\documentclass{amsart}
\usepackage{amssymb}
\usepackage{amsfonts}
\usepackage{amsmath}
\usepackage{lastpage}
\usepackage[legalpaper,bookmarks=true,colorlinks=true,linkcolor=blue,citecolor=blue]{hyperref}
\usepackage{graphicx}%
\usepackage{fancyhdr}
\usepackage{color}
\usepackage[mathlines]{lineno}
\usepackage{lscape}
\usepackage{epsfig}
\usepackage{natbib}
\usepackage{geometry}
\usepackage{tgbonum}
\fontfamily{qcr}\selectfont

\newtheorem{theorem}{Theorem}
\theoremstyle{plain}

\newtheorem{lemma}{Lemma}

\newtheorem{proposition}{Proposition}

\numberwithin{equation}{section}

\newcommand{\Bin}{\bigskip \noindent}

\newcommand{\Ni}{\noindent}

\begin{document}
\Large
\title{A detailed and direct proof of Skorohod-Wichura's theorem}
\author{Gane Samb Lo $\dag$}
\author{$\dag\dag$ Aladji Babacar Niang}
\author{Lois Chinwendu Okereke $\dag\dag\dag$}

\begin{abstract} The representation Skorohod theorem of weak convergence of random variables on a metric space goes back to Skorohod (1956) in the case where the metric space is the class of real-valued functions defined on [0,1] which are right-continuous and have left-hand limits when endowed with the Skorohod metric. Among the extensions of that to metric spaces, the version by Wichura (1970) seems to be the most fundamental. But the proof of Wichura seems to be destined to a very restricted public. We propose a more detailed proof to make it more accessible at the graduate level. However we do far more by simplifying it since important steps in the original proof are dropped, which leads to a direct proof that we hope to be more understandable to a larger spectrum of readers. The current version is more appropriate for different kinds of generalizations.

\bigskip
\noindent $^{\dag}$ Gane Samb Lo.\\
LERSTAD, Gaston Berger University, Saint-Louis, S\'en\'egal (main affiliation).\newline
LSTA, Pierre and Marie Curie University, Paris VI, France.\newline
AUST - African University of Sciences and Technology, Abuja, Nigeria\\
gane-samb.lo@edu.ugb.sn, gslo@aust.edu.ng, ganesamblo@ganesamblo.net\\
Permanent address : 1178 Evanston Dr NW T3P 0J9,Calgary, Alberta, Canada.\\

\noindent $^{\dag\dag}$ Aladji Babacar Niang\\
LERSTAD, Gaston Berger University, Saint-Louis, S\'en\'egal.\\
Email: aladjibacar93@gmail.com\\

\noindent $^{\dag\dag\dag}$ Lois Chinwendu Okereke\\
Department of Pure and applied Mathematics\\
African University of Sciences and Technologies (Abuja), Abuja, Nigeria\\
Email : lokereke@aust.edu.ng\\

%\bigskip
%\noindent $^{\dag\dag}$ Mohammad Ahsanullah\\
%Department of Management Sciences. Rider University. Lawrenceville, New Jersey, USA\\
%Email : ahsan@rider.edu\\

%\noindent \noindent $^{\dag\dag\dag}$ Modou Ngom.\\
%LERSTAD, Gaston Berger University, Saint-Louis, S\'en\'egal (main affiliation).\newline
%modou.ngo.edu@gmail.com\\

%\noindent \noindent $^{\dag\dag\dag}$ Modou Ngom.\\
%LERSTAD, Gaston Berger University, Saint-Louis, S\'en\'egal (main affiliation).\newline
%modou.ngo.edu@gmail.com\\

\end{abstract}
\maketitle

\section{introduction}

\subsection{An easy entry} Let us begin by the simplest form of the Skorohod theorem which is the object of our study. By the Portmanteau Theorem (see \cite{ips-wcia-ang}, Corollary 3, page 72), the weak convergence of probability measures on $\mathbb{R}$ is reduced to the convergence of cumulative distribution functions (\textit{cdf}). Precisely, the sequence of probability measures $(\mathbb{P}_n)_{n\geq 1}$ on the Borel space $(\mathbb{R}, \mathcal{B}(\mathbb{R}))$ associated with the \textsl{cdf}'s

$$
F_n(x)=\mathbb{P}_n(]-\infty,x]), \ x \in \mathbb{R}, \ \ n \geq 1,
$$

\Bin weakly converges to the probability measure $\mathbb{P}_{\infty}$ associated with the \textsl{cdf} $F_{\infty}(x)=\mathbb{P}_{\infty}(]-\infty,x]), \ x \in \mathbb{R}$, denoted by $\mathbb{P}_n \rightsquigarrow \mathbb{P}_{\infty}$, in the sense that for all $f : \mathbb{R} \rightarrow \mathbb{R}$ continuous and bounded

\begin{equation}
\int f \ d\mathbb{P}_n \rightarrow  \int f \ d\mathbb{P}_{\infty}, \ as \ n\rightarrow +\infty \label{Def1}
\end{equation}

\noindent if and only if, for any continuity point $x$ of $F_{\infty}$ (denoted $x \in C(F_{\infty}))$,

\begin{equation}
F_n(x) \rightarrow F_{\infty}(x), \ as \ n \rightarrow +\infty.  \label{def2}
\end{equation}

\bigskip \noindent The convergence in Formula \eqref{def2} is still denoted by $F_n \rightsquigarrow F_{\infty}$. That formula itself becomes a definition of weak convergence and the notion is generalized to convergence of sequences of non-decreasing functions to another non-decreasing function. By using the generalized inverse  (see \cite{ips-wcia-ang}, page 112) of a non-decreasing function $F$ defined by

\begin{equation}
F^{-1}(y)= \inf \{x \in \mathbb{R}, \ F(x) \geq y\}, \ y \in \mathbb{R}, 
\end{equation}

\bigskip 
\noindent which is non-decreasing too (and is the inverse function of $F$ if $F$ is invertible), \cite{billingsley1968} proved that

\begin{equation}
F_n \rightsquigarrow F_{\infty} \Rightarrow F_n^{-1} \rightsquigarrow F_{\infty}^{-1}. \label{def3}
\end{equation}

\bigskip 
\noindent When applied to random variables $(X_n)_{n\geq 1}$ and $X_{\infty}$ taking values in $\mathbb{R}$, we consider the probability laws 
and the \textit{cdf}'s of those probability laws denoted by $F_{X_{\infty}}$ and $(F_{X_n})_{n\geq 1}$. Then Formula \eqref{def2} becomes

\begin{equation}
\forall x \in C(F_{X_\infty}), \ F_{X_n}(x) \rightarrow F_{X_\infty}(x), \label{def4}
\end{equation}

\bigskip \noindent and we still write $X_n \rightsquigarrow X_{\infty}$.\\

\noindent The objective of the Skorohod problem is the following. In Formula \eqref{def4}, only the probability laws of $X_{\infty}$ and the $X_n$'s matter. So the probability space on which $X_\infty$ or each $X_n$, $n\geq 1$, is defined has no importance for the convergence. They may even be pair-wise distinct. As a result, we neither have access to the paths of the random variables, nor can we make coherent numerical operations on them. On the other hand, it is also known that if $X_\infty$ and all the  $X_n$'s are on the same probability space, the almost-sure convergence and the convergence in probability imply the weak convergence. The following question arises :\\

\noindent (\textbf{Q}) Given \eqref{def4}, can we construct a probability space $(\Omega, \mathcal{A}, \mathbb{P})$ holding a modification $\{Y_{\infty}, \ Y_n, \ n\geq 1\}$
of $\{X_{\infty}, \ X_n, \ n\geq 1\}$, in the sense that each $Y_n$ has the same law as $X_n$ for $n \in \mathbb{N} \cup \{+\infty\}$, such that the weak convergence of $Y_n$ to $Y_{\infty}$ is an almost-sure convergence?\\

\noindent A successful answer to this question allows to treat problems of weak convergence as almost-sure convergences and  to use analysis tools such as expansions, delta-methods, etc., in solving weak convergence problems.\\

\noindent The answer is relatively easy on $\mathbb{R}$ due to Formula \eqref{def3} based on Renyi's representations of real valued random variables
(see \cite{ips-wcia-ang}, page 127). Indeed we may take $(\Omega, \mathcal{A}, \mathbb{P})=([0,1], \mathcal{B}([0,1]), \lambda)$, where $\lambda$ is the Lebesgue measure which actually is a probability measure on $[0,1]$. Let $U : [0,1] \rightarrow [0,1]$ be the identity mapping. It is a real-valued random variable defined on $[0,1]$ which is uniformly distributed on $(0,1)$. By the Renyi representation,  $\{F_{X_\infty}^{-1}(U), \ F_{X_n}^{-1}(U), \ n\geq 1\}$ is a modification of $\{X_{\infty}, \ X_n, \ n\geq 1\}$. Using \eqref{def3}, we have

\begin{eqnarray*}
&&\mathbb{P}(\{\omega \in \Omega, F_{X_n}^{-1}(U(\omega)) \nrightarrow F_{X_\infty}^{-1}(U(\omega))\}\\
&&=\lambda(\{u \in [0,1],  F_{X_n}^{-1}(u) \nrightarrow F_{X_\infty}^{-1}(u))\}\\
&&\leq \lambda(C(F_\infty^{-1})^c)=0,
\end{eqnarray*}

\bigskip 
\noindent since $F_\infty^{-1}$ is non-decreasing and hence has at most a countable number of discontinuity points and countable sets are null sets for the Lebesgue measure. $\noindent$\\

\noindent This is the Skorohod theorem on $\mathbb{R}$.\\

\subsection{More general versions}

\noindent The Skorohod theorem as well as the Skorohod topologies go back to the paper by \cite{skorohod1956} in which the theorem is stated on the class $D([0,1])$ of real-valued functions defined on $[0,1]$ which are right-continuous and have left-hand limits, endowed with the topology he created and that is named after him. Since then, many authors have tried to extend it to general metric spaces, at least separable. (see \cite{dudley66a} for non-separable spaces).\\

\noindent The version given by \cite{wichura1970} can be considered as of the most fundamental, from which further versions can be done. The proof of \cite{wichura1970} is indirect and certainly reserved to a very restricted public.\\

\noindent Given the importance of that theorem in weak convergence, we think that it deserves to have a more direct and pedagogical version. We want to make the proof available to readers coming first from courses of Topology, Measure and Integration, Mathematical Probability Theory but mastering the product of measurable spaces and product of measures, the beginners of graduate studies in summary. This what we intend to do here.\\

\noindent Before we proceed, we advise readers not familiar with arbitrary products of probability spaces, complete probability spaces and/or Caratheordory extensions of measures on semi-algebra to read about them in appropriate books, possibly in Chapter 9 in \cite{ips-mfpt-ang}, Chapter 5 in \cite{ips-mestuto-ang} (see part II on Outer measures and Doc. 04-02, Exercise 8).\\

\noindent Here, we mainly follow the ideas in the proofs of \cite{wichura1970}, which is valid in any metric space, when the probability measure limit is tight.  However, our modifications are significantly interesting in terms of new results concerning a number of assertions. We tried to clarify many details at a cost of a  few more pages. The organization of the proofs is radically different. We will comment on the key points of our methods compared to the initial proof.\\

\noindent Here is the general statement of the main theorem.\\

\begin{theorem} \label{wcib_SkororhodWhichura}(Skorohod-Wichura) Let $D$ be an upward well-directed set, $(\mathbb{P}_{\alpha})_{\alpha \in D}$ be a family of probability laws on a metric space $(S,d)$ endowed with its Borel $\sigma$-algebra $\mathcal{S}=\mathcal{B}(S)$ and let $\mathbb{P}_{\infty}$ be a tight probability law on $(S,d)$. Suppose that we have

\begin{equation}
\mathbb{P}_{\alpha} \rightsquigarrow \mathbb{P}_{\infty} \ on \ D. \label{wcib_weakconv}
\end{equation}

\bigskip \noindent Then, there exists a probability space $(\Omega, \mathcal{A}, \nu)$ holding a family of random variables $(X_{\alpha})_{\alpha \in D}$ and a random variable $X_{\infty}$ taking values in $S$ such that :\\

\noindent (a) $\nu_{X_{\infty}}=\nu X_{\infty}^{-1}=\mathbb{P}_{\infty}$,\\

\noindent (b) $\forall \alpha \in D$, $\nu_{X_{\alpha}}=\nu X_{\alpha}^{-1}=\mathbb{P}_{\alpha}$\\

\noindent and\\ 

\noindent (c) the field $(X_{\alpha})_{\alpha \in D}$ almost surely converges to $X_{\infty}$.
\end{theorem}

\bigskip

\noindent We will need some preliminary tools and remarks.

\section{Preliminary tools}  \label{07_sec_02_ss1}

\noindent By assumption, $\mathbb{P}_{\infty}$ is tight. So for any integer $n\geq 1$, there is a compact subset $K_n$ of $S$ such that $\mathbb{P}_{\infty}(K_n)\geq 1-1/n$. By taking $S_{\infty}$ as the closure of $\cup_{n\geq 1} K_n$, we have

$$
\mathbb{P}_{\infty}(S_{\infty})=1.
$$

\Bin This means that $S_{\infty}$ a $\sigma$-compact set (a countable union of compacts sets) and that $S_{\infty}$ is a separable. We may replace the weak limit $\mathbb{P}_{\alpha} \rightsquigarrow \mathbb{P}_{\infty}$ by $\mathbb{P}_{\alpha}(\circ \cap S_{\infty})\rightsquigarrow \mathbb{P}_{\infty}(\circ \cap S_{\infty})$, meaning that the convergence holds on $S_{\infty}$.\\

\noindent So we can do proceed to the proof when $S$ is a separable and complete set (Polish space) on which each probability measure is tight. To begin, we give the following first tool.\\

\begin{proposition} \label{07_skorohod_propPolish} Let $\mathbb{P}_{\infty}$ be a probability measure on the separable space $(S,d)$. Let $\Delta>0$ and $\varepsilon>0$. Then for any $k\geq 1$, there exists a finite partition $C_{j,k}$, $0\leq j \leq q(k)$ of $\mathbb{P}_{\infty}$-continuous borel set such that

$$
\forall (k\geq 1), \ \mathbb{P}_{\infty}(C_{0,k})\leq 2^{-k}\varepsilon=:\varepsilon_k $$

\Bin and

$$
\forall (k\geq 1), (\forall 1\leq j \leq q(k)),  \ diam(C_{j,k})\leq 2^{-k}\Delta.
$$

\bigskip \noindent Moreover, for all $k\geq 1$, any element of the partition $C_{j,k+1}$, $0\leq j \leq q(k+1)$ is the sum of $\mathbb{P}_{\infty}$-continuous sets, each of them being included in one element of the partition $C_{j,k}$, $0\leq j \leq q(k)$. 
\end{proposition}

\noindent \textbf{Proof of Proposition \ref{07_skorohod_propPolish}}. Since $S$ is separable, by the Lindel\H{o}f property, any open cover can be reduced to one of its countable sub-covers (See \cite{choquet1966}). Let $D$ be a dense subset of $S$. For a fixed real number $\Delta>0$, we may cover $S$ by a countable number of balls centered on points forming a subset $D_0$ of $D$ with radius $\Delta/2$. For each  $s \in D_0$, the borders of the balls $B(s,\delta +\Delta/4)$, $0\leq \delta \leq \Delta/4$ are disjoint. Thus, we can find a value $\delta$ such that $\mathbb{P}_{\infty}(\partial B(s,\delta +\Delta/2))=0$ (See \cite{ips-mestuto-ang} ). So 
$S$ may be covered by a countable $\mathbb{P}_{\infty}$-continuous balls of diameters at most equal to $\Delta_1=\Delta$. If we denote these balls by $D_{1,n}$, $n\geq 1$, we have

$$
S=\bigcup_{n\geq 1} D_{1,n}.
$$ 

\bigskip \noindent By the continuity of the probability measure, there exists $n(1)\geq 1$ such that

$$
\mathbb{P}_{\infty}\biggr(\bigcup_{1 \leq n \leq n(1)} D_{1,n} \biggr)\geq 1-\varepsilon_1.
$$

\bigskip \noindent Now we can use the classical method to transform the union of the  $D_{1,n}$, $1\leq n \leq n(1)$ into a sum of sets, that is

$$
S_1=\bigcup_{1 \leq n\leq n(1)} D_{1,n}=\sum_{1 \leq n\leq n(1)} C_{n,1}, 
$$

\bigskip \noindent with

$$
C_{1,1}=D_{1,1}, \ C_{2,1}=D_{1,1}^{c}D_{1,2} \subset D_{1,2}, \ C_{j,1}=D_{1,1}^{c}...D_{1,j-1}^{c}D_{1,j}\subset D_{1,j}, \ \ j\geq 3,
$$

\bigskip \noindent (with $q(1)=n(1))$. By the following properties  [where int(A) and ext(A) stand for the interior and of the exterior of a set $A$, respectively], 

$$
\partial (A)= (int{A} \cap ext(A))^c, \ int{A^c}=ext(A) \ \ and \ \ ext(A^c)=int(A), 
$$

\bigskip \noindent we can see that a set and its complement has the same border and we already know that the border of a countable union of sets is included in the union of the borders of the sets. When combined, the two points  say that intersections and unions of sets are included in the union of the borders of those sets. So all the $C_{n,1}$, $1\leq n \leq n(1)$, are disjoint and $\mathbb{P}_{\infty}$-continuous measurable sets of diameters at most equal to $\Delta_1=\Delta$. By putting

$$
C_{0,1}=S \setminus \bigcup_{1 \leq n\leq n(1)} D_{1,n},
$$

\bigskip \noindent we still get a $\mathbb{P}_{\infty}$-continuous measurable set. So we get that $S$ is partitioned into the disjoint measurable sets $C_{n,1}$, $0 \leq n \leq n(1)$, all $\mathbb{P}_{\infty}$-continuous and having diameters  at most equal to $\Delta_1=\Delta$ such that 

$$
\mathbb{P}_{\infty}(C_{0,1}) \leq \varepsilon_1.
$$

\bigskip \noindent In a second step, let us say that each $D_{1,j}$, $0 \leq j \leq n(1)$, may be covered (in $S$) by a countable of balls centered on elements of some $D_2(j) \subset D \cap D_{1,j}$ (we do not forgot that $D_{1,j}$ is an open ball) with radius at most $\Delta/4$. We repeat the argument in the first step by using the balls $B(s,\delta+\Delta/8)$, $s\in D_2(j)$, $0\leq \delta \leq \Delta/8$ and choose $\delta_{s}$ such that $B(s,\delta_s+\Delta/8)$ is $\mathbb{P}_{\infty}$-continuous and we get

$$
D_{1,j}=D_{1,j} \bigcap \bigcup_{s \in D_2(j)} B(s,\delta_s+\Delta/4)\equiv \bigcup_{s\in D_2(j)} D_{s,j}, 
$$  

\bigskip \noindent where the  $D_{s,j}$'s are $\mathbb{P}_{\infty}$-continuous of diameter at most equal to $\Delta/2$. Se have 

$$
S_1=\bigcup_{1\leq j \leq n(1)} \bigcup \bigcup_{s \in D_2(j)} D_{s,j} =: \bigcup_{n\geq 1} D_{2,n},
$$

\bigskip \noindent where the $D_{2,n}$ are open sets, all $\mathbb{P}_{\infty}$-continuous of diameter at most equal to $\Delta/2$. From there, we proceed as in the conclusion of the first step to get measurable sets  $C_{n,2}$, $1\leq n \leq n(2)=q(2)$, all $\mathbb{P}_{\infty}$-continuous of diameter at most equal to $\Delta/2$ such that for 

$$
C_{0,2} =S_1 \setminus \left(\bigcup_{1\leq n \leq q(2)} D_{2,n}\right),
$$

\Bin we have

$$
\mathbb{P}_{\infty}(C_{0,2})\leq \varepsilon_2,
$$

\Bin and surely, for each $j \in \{1,\cdots,q(2)\}$, $C_{j,2}$ is in one of the $C_{h,1}$'s.\\

\bigskip \noindent From this the proof of the proposition is completed by induction. $\blacksquare$\\

\Ni We will also need the following result.\\

\begin{proposition} \label{wcib_07_countable_Product}
Let $(\Omega_D=\prod_{\alpha \in D} S_{\alpha},\ \mathcal{S}_{D}=\otimes_{\alpha \in D} \mathcal{S}_{\alpha}, \nu)$ be a non-countable product space holding a $\sigma$-finite measure which takes finite values on each element of a countable partition of $ \prod_{\alpha \in D} S_\alpha$ consisting of cylinders. Let 
$A \in \otimes_{\alpha \in D} \mathcal{S}_{\alpha}$ be a measurable set.\\

\noindent Then $A$ is a subset of a measurable set $B$ depending only on a countable number of factors $\alpha \in D$. 
\end{proposition}

\bigskip \noindent \textbf{Proof of Proposition \ref{wcib_07_countable_Product}}. It is enough to do the proof with a finite measure $\nu$ since the extension to a $\sigma$-finite measure is straightforward. Hence $\nu$ is a measure on the algebra $\mathcal{C}_D$ of finite sums of cylinders. A measurable set $A$ is of the form

$$
A=\prod_{\alpha \in V} A_{\alpha} \times \prod_{\alpha \notin V} S_{\alpha} \equiv A_{V} \times S^{\prime}_{V},  \ \ \ (P1)
$$ 

\bigskip
\noindent where $V$ is a finite subset of $D$ and $A_\alpha \in \mathcal{S}_{\alpha}$. We already suggested the reader to read Chapter 8 in \cite{ips-mfpt-ang} about product $\sigma$-algebras and relevant notation. Now, cylinders, as denoted in Formula (P1), depends only on a finite number of indices, since  $x=(x_{\alpha})_{\alpha \in D}$ is in $A$ if and only if $x_V=(x_{\alpha})_{\alpha \in V}$ is in $A_V$.\\

\noindent Let us recall the projections on sub-products space of $\Omega_D$. Define for any subset $D_1\subset D$, the projection $\pi_{D_1}$ of 
$\prod_{\alpha \in D} S_{\alpha}$ on $\prod_{\alpha \in D_1} S_{\alpha}$, defined by 

$$
\forall x=(x_\alpha)_{\alpha \in D} \in \prod_{\alpha \in D} S_{\alpha}, \ \pi_{D_1}(x)=(x_\alpha)_{\alpha \in D_1}.
$$

\Bin  Let us proceed to the extension of the measure $\nu$ from the algebra $\mathcal{C}_D$ to the measure $\tilde{\nu}$ on $\sigma$-algebra $\mathcal{S}_D$ it generates  by the outer-measure method. But $\tilde{\nu}$ and $\nu$ are equal on $\mathcal{S}_{D}$. So by the exterior measure definition (See \cite{ips-mestuto-ang}, Chapter 5, Doc 04-10), we have

$$
\nu(A)=\inf\left\{\nu(B), A \subset B=\bigcup_{n\geq 1} B_n, \ B_n \in \mathcal{C}_D\right\}.
$$

\bigskip \noindent There is nothing to do when $A=\Omega_D$. If $\nu(A)<\nu(\Omega)$, for $\varepsilon<(\nu(\Omega)-\nu(A))/2$, there exists a countable number of finite sums of cylinders, denoted by $B_n$, $n\geq 1$ such that, for $B$ being the unions of the $B_n$'s,

$$
A \subset B \ \ \ and  \ \ \ \nu(A)\leq  \nu(B) < \nu(A)+\varepsilon.
$$ 

\bigskip \noindent To conclude, we take as $D_0$ the set of indices $\alpha$ involved in the cylinders forming the finite sums of cylinders $B_n$, $n\geq 1$. We have that $B$ is of the form

$$
B=B_{D_0} \times S^{\prime}_{D_0}. \ \ \ \square.
$$

\bigskip \noindent Now we give the complete proof of the main theorem.\\

\Ni \textbf{Proof of Theorem \ref{wcib_SkororhodWhichura}}.  \label{07_sec_02_ss2}  Now, let us describe the space on which holds our construction. We replicate the measurable space $(S, \mathcal{S})$ into $(S_{\infty}, \mathcal{S}_{\infty})$, $(S_{\alpha}, \mathcal{S}_{\alpha})$, $\alpha \in D$, in sense that each $(S_{\alpha}, \mathcal{S}_{\alpha})$ is identical to $(S, \mathcal{S})$ for $\alpha \in \{\infty\} \cup D$. We consider the measurable product space

$$
\Omega =\prod_{\alpha \in \{\infty\} \cup D} \ \ S_{\alpha}=S_{\infty} \times \prod_{\alpha \in  D} \ S_{\alpha} 
$$

\Bin and endow it with the product $\sigma$-algebra

$$
\mathcal{A}=\otimes_{\alpha \in \{\infty\} \cup D} \ \mathcal{S}_{\alpha}= \mathcal{S}_{\infty} \otimes \otimes_{\alpha \in D} \ \mathcal{S}_{\alpha}.
$$

\Bin We take $(\beta_k)_{k\geq 1}$ an arbitrary sequence such that each $\beta_k>0$, 

$$
\beta_1>\beta_2>\cdots > \beta_k \downarrow 0 \ as \ k\uparrow +\infty \ \ and \ \ \sum_{k\geq 1}\beta_k=1,
$$

\Bin (we may take   $\beta_k=2^{-k}$, $k\geq 1$, for example) and we denote $\beta_{k}^{\star}=\beta_{1}+\cdots+\beta_{k}$, $k\geq 1$. We denote elements of $\Omega$ as

$$
\omega=(\omega_{\infty},  (\omega_\alpha)_{\alpha \in D}).
$$

\Bin The simple projections $\pi_{\infty}$ and $\pi_{\gamma}$, for $\gamma \in D$, are defined as

$$
\pi_{\infty}(\omega)=\pi_{\infty}(\omega_{\infty}, (\omega_\alpha)_{\alpha \in D})=\omega_{\infty}.
$$

\noindent and

$$
\pi_{\gamma}(\omega)=\pi_{\gamma}(\omega_{\infty}, (\omega_\alpha)_{\alpha \in D})=\omega_{\gamma}
$$

\Bin We adopt the usual notation like in \cite{ips-mfpt-ang} (Chapter 9). For example, for a non-empty subset $D_1 \subset \{\infty\} \cup D$, for 

$$
A_{D_1} \subset \prod_{\alpha \in D_1} S_\alpha
$$  

\Bin the cylinder of base $A_{D_1}$ is denoted as

$$
A_{D_1} \times \prod_{\alpha \notin D_1} \ S_{\alpha}=\{(\omega_{\infty}, (\omega_\alpha)_{\alpha \in D}), \ (\omega_\alpha)_{\alpha \in D_1} \in A_{D_1}\},
$$

\Bin and we denote

$$
A_{D_1} \times \prod_{\alpha \notin D_1} \ S_{\alpha}=A_{D_1} \times S^{\prime}_{D_1}.
$$

\Bin When we consider the measurable space $(\Omega, \mathcal{A})$, the identity mapping

$$
X \ : \ (\Omega, \mathcal{A}) \rightarrow \biggr(S_{\infty} \times \prod_{\alpha \in  D} \ S_{\alpha}, \  \mathcal{S}_{\infty} \otimes \otimes_{\alpha \in D} \ \mathcal{S}_{\alpha}\biggr)
$$

\Bin is measurable and by the way, for any $\omega=(\omega_{\infty},  (\omega_\alpha)_{\alpha \in D}) \in \Omega$,

$$
X(\omega)=\left(X_{\alpha}(\omega)\right)_{\alpha \in \{\infty\} \cup D}.
$$

\Bin In the current case of the identity function, we have

$$
X_{\infty}(\omega)=\omega_{\infty} \ \ and \ \ X_{\alpha}(\omega)=\omega_{\alpha} \ for \ \alpha \in D
$$

\Bin and for $A_{\infty} \subset S_{\infty}$ and $A_{\alpha} \subset S_{\alpha}$ for $\alpha \in D$, we have by the same notation,

$$
X_{\infty}^{-1}(A_{\infty})=A_{\infty} \times S^{\prime}_{\{\infty\}} \ and \ X_{\alpha}^{-1}(A_{\alpha})=A_{\alpha} \times S^{\prime}_{\{\alpha\}}
$$

\Bin Later, we will create a probability measure $\nu$ on $(\Omega, \mathcal{A})$ such that $\nu X_{\alpha}^{-1}=\mathbb{P}_{\alpha}$, 
$\alpha \in \{\infty\} \cup S$.\\

\Ni After this important notation, we begin by exploiting the partitions and facts in Proposition \ref{07_skorohod_propPolish} on each $S_{\alpha}=S$ (for $\alpha \in \{\infty\} \cup D$) and we denote $\Delta_k=\Delta/2^k$, $\varepsilon_k=2^{-k} \varepsilon$. For each $k\geq 1$, we denote

$$
I(k)=\{0\leq j \leq q(k), \ \mathbb{P}_{\infty}(C_{j,k})\neq 0\}.
$$ 

%We have $0 \in I(k)$ since $\mathbb{P}_{\infty}(C_{0,k})\geq 1-\varepsilon_k$.
\bigskip \noindent  By the weak convergence of $\mathbb{P}_{\alpha}$ to $\mathbb{P}_{\infty}$, and by the $\mathbb{P}_{\infty}$-continuity of the $C_{j,k}$'s, we have for all  $j \in \{0,\cdots,q(k)\}$, that

$$
\mathbb{P}_{\alpha}(C_{j,k}) \rightarrow \mathbb{P}_{\infty}(C_{j,k}).
$$

\bigskip \noindent So, for a fixed $c \in ]0,1[$, for any $k\geq 1$,  there exists $\alpha_k$ such that for $\alpha\geq \alpha_k$, for all $j\in I(k)$, 

$$
\mathbb{P}_{\infty}(C_{j,k})>0 \ \ and \ \ 1-c \leq \frac{\mathbb{P}_{\alpha}(C_{j,k})}{\mathbb{P}_{\infty}(C_{j,k})}\leq 1+c. \ \ (G1)
$$

\bigskip \noindent We may suppose (G1) holds for all $\alpha \in D$ for $k=1$. Now let us define

\begin{equation}
\mathbb{H}_{\alpha,k}(\circ)=\frac{1}{1-\beta_{k}^{\ast}}\mathbb{P}_{\alpha}(\circ)- \sum_{j\in I(k)}\frac{\beta_{k}^{\ast}}{1-\beta_{k}^{\ast}}\mathbb{P}_{\alpha}(\circ/C_{j,k})\mathbb{P}_{\infty}(C_{j,k}). \label{wcib_defHalpha}
\end{equation}

\bigskip \noindent It is clear that each $\mathbb{H}_{\alpha,k}$ is $\sigma$-additive and $\mathbb{H}_{\alpha,k}(S_{\alpha})=1$. Further, for any $A_{\alpha} \in \mathcal{S}_{\alpha}$, we have

$$
\mathbb{H}_{\alpha,k}(A_{\alpha})=\sum_{p\in I(k)} \mathbb{H}_{\alpha,k}(A_{\alpha}\cap C_{p,k}),
$$

\bigskip \noindent and for all $p\in I(k)$,

\begin{eqnarray*}
\mathbb{H}_{\alpha,k}(A_{\alpha}\cap C_{p,k})&=&\frac{1}{1-\beta_{k}^{\ast}}\mathbb{P}_{\alpha}(A_{\alpha}\cap C_{p,k})- \frac{\beta_{k}^{\ast}}{1-\beta_{k}^{\ast}}
\mathbb{P}_{\alpha}(A_{\alpha}/C_{p,k})\mathbb{P}_{\infty}(C_{p,k})\\
&=&\frac{1}{1-\beta_{k}^{\ast}}\mathbb{P}_{\alpha}(A_{\alpha} \cap C_{p,k})- \frac{\beta_{k}^{\ast}}{1-\beta_{k}^{\ast}}
\mathbb{P}_{\alpha}(A_{\alpha}\cap C_{p,k}) \frac{\mathbb{P}_{\infty}(C_{p,k})}{\mathbb{P}_{\alpha}(C_{p,k})}.
\end{eqnarray*}

\noindent Thus

\begin{equation*}
\mathbb{H}_{\alpha,k}(A_{\alpha}\cap C_{p,k})=\mathbb{P}_{\alpha}(A_{\alpha}\cap C_{p,k})\biggr(\frac{1}{1-\beta_{k}^{\ast}}-\frac{\beta_{k}^{\ast}}{1-\beta_{k}^{\ast}}\frac{\mathbb{P}_{\infty}(C_{p,k})}{\mathbb{P}_{\alpha}(C_{p,k})}\biggr),
\end{equation*}

\bigskip \noindent and the term between the parentheses is non-negative if and only if

$$
\beta_{k}^{\ast} \leq \frac{\mathbb{P}_{\alpha}(C_{p,k})}{\mathbb{P}_{\infty}(C_{p,k})}\equiv \eta_{\alpha, k,p}. \ \ \ (G2)
$$

\bigskip \noindent Hence $\mathbb{H}_{\alpha,k}$ is non-negative whenever Formula (G2) holds for all $p\in I(k)$. From there, we re-use the method of the first part. We consider sequences

$$
1-c=\varepsilon_{1,0} < \varepsilon_{1,1}<\cdots < \varepsilon_{1,r} \uparrow 1, \ as \ r\uparrow +\infty
$$

\noindent and

$$
\varepsilon_{2,r} <\varepsilon_{2,r-1} <  \cdots <\varepsilon_{2,1}< \varepsilon_{2,0}=1+c,  \ \varepsilon_{2,r}  \downarrow 1, \ as \ r\uparrow +\infty.
$$

\bigskip \noindent  For any fixed $k\geq 1$, we define a mapping $\ell_1 : \ D \rightarrow \mathbb{N}^\ast\cup \{+\infty\}$ case by case.\\

\noindent Case 1 : Let us denote $\eta_{\alpha, k}=\min_{p \in I(k)} \eta_{\alpha, k,p}$. For $\alpha \in D$ fixed, if all the $\eta_{\alpha, k}$'s, are less that $1-c$, we put  $\ell_1(\alpha)=+\infty$.\\

\noindent Case 2 : if all the $\eta_{\alpha, k}$'s are greater that $1+c$, we take $\ell_1(\alpha)=+\infty$.\\

\noindent Case 3 : if one of the two cases above fails, we may define

$$
J_1(\alpha)=\{j\geq J, \ \exists (k\geq 1), \  \eta_{\alpha, k} \in ]\varepsilon_{1,j}, \varepsilon_{1,j+1}]\}
$$

\Bin and

$$
J_2(\alpha)=\{j\geq 1, \ \exists (k\geq 1), \ \eta_{\alpha, k} \in ]\varepsilon_{2,j+1}, \varepsilon_{2,j}]\}.
$$

\Bin For $i \in \{1,2\}$, if $J_i(\alpha)$ is empty or non-empty but unbounded, we put $\hat \ell(i,\alpha)=+\infty$ and if it is non-empty and bounded we set
$j_i(\alpha)=\max J_i(\alpha)$ (for $\tau_i=1_{(i=2)})$

$$
\hat  \ell(i,\alpha)=\underset{\beta_{k}^{\ast} \leq \varepsilon_{i,j_i(\alpha)+\tau_j}}{Argmax \ } \beta_{k}^\ast.
$$

\Bin In all cases, we take

$$
\ell_1(\alpha)=min(\hat  \ell(1,\alpha), \ \hat  \ell(2,\alpha))
$$

\Bin Conclusion, the function $\ell_1(\circ)$ is well-defined. Let us construct a second mapping $\ell_2$ as follows. We know that for any fixed $k\geq 1$, we have

$$
\delta_{\alpha,k}=\max_{0\leq h \leq q(k)} |\mathbb{P}_{\alpha}(C_{h,k})-\mathbb{P}_{\infty}(C_{h,k})|\rightarrow 0 \ on \ D.
$$

\Bin From the the ranking $\beta_1>\beta_2>\cdots>\beta_j \searrow 0$, we define the non-empty set

$$
J_3(\alpha)=\{ j\geq 1, \ (\exists k\geq 1), \ \delta_{\alpha,k} \in ]\beta_{j+1}, \beta_{j}]\}
$$

\Bin and take $\ell_2(\alpha)=+\infty$ if $J_3(\alpha)$ is unbounded and $\ell_2(\alpha)=\max J_3(\alpha)$ otherwise. Finally, we take

$$
\ell=\min(\ell_1, \ell_2).
$$

\bigskip \noindent Let us show that $\ell(\alpha) \rightarrow +\infty$ on $D$. Let $k_0>0$, and consider the unique value $j_1$ of $j$ such that  $\beta_{k_0}^\ast \in ]\varepsilon_{1,j}, \varepsilon_{1,j+1}]$. All the $\eta_{\alpha, k_0,p}$'s converge to one and all the $\delta_{\alpha,k_0}$'s converge to zero. So, for some $\alpha_{\infty}$, for all $\alpha\geq \alpha_{\infty}$, we have that $\eta_{\alpha, k_0}$ is in some interval $]\varepsilon_{1,r},\varepsilon_{2,r}]$, where $r\geq j_1$, and all the $\delta_{\alpha,k_0}$'s are in some interval $]\beta_{s+1}, \beta_{s}]$, $s\geq k_0$. By definition of $\ell_1$, in Case 3, and by definition of $\ell_2$, we get $\ell(\alpha)\geq \min(r,s)\geq k_0$. We have :

$$
\exists \ \ell : \ D \rightarrow \mathbb{N}^\ast\cup \{+\infty\}, \ \lim_{D} \ell(\alpha)=+\infty.
$$

\Bin Moreover, by construction, we have that, for all $\alpha \in D$ such that $\ell(\alpha)<+\infty$,

$$
\mathbb{H}_{\alpha}=:\mathbb{H}_{\alpha,\ell(\alpha)},
$$

\bigskip \noindent is non-negative and hence is a probability measure such that (by inverting the formula)

\begin{equation} \label{wcib_defHalphaB}
\mathbb{P}_{\alpha}(\circ)=(1-\beta_{\ell(\alpha)}^\ast) \mathbb{H}_{\alpha}(\circ)+\beta_{\ell(\alpha)}^\ast \sum_{p \in I(\ell(\alpha))} \mathbb{P}_{\infty}(C_{p,\ell(\alpha)}) \mathbb{P}_{\alpha}\left(\circ/C_{p,\ell(\alpha)}\right)
\end{equation}

\noindent holds and we also have \label{AP0}

$$
|\mathbb{P}_{\alpha}(C_{h,\ell(\alpha)})-\mathbb{P}_{\infty}(C_{h,\ell(\alpha)})|\geq \beta_{\ell(\alpha)}. \ \ \ \  (AP0)
$$

\bigskip \noindent Let us point out that $\mathbb{H}_\alpha$ is used in all this chapter only if $\ell(\alpha)<\infty$. Now we may construct the space on which the almost-sure limit will holds.\\

\Bin For any $\alpha \in D$, for $s \in S$,  by the decomposition of $S$ into the $C_{p,\ell(\alpha)}$'s, $0\leq p\leq q(\ell(\alpha))$, there exists a unique $p(s,\alpha)$ such that $s \in C_{p(s,\alpha),\ell(\alpha)}$. We define for $\alpha\in D$, $1\leq j <+\infty$, $0\leq p \leq q(\ell(\alpha))$,

\begin{equation*}
\nu_{j,s,\alpha}=\left\{ 
\begin{tabular}{lll}
$\mathbb{H}_{\alpha}$ & if & $j>\ell(\alpha)$ \\ 
$\mathbb{P}_{\alpha}(\circ /C_{p(s,\alpha),\ell(\alpha)})$ & if  & $ j\leq \ell(\alpha)$.\\
\end{tabular}
\right. ,
\end{equation*}

\Bin for $\ell(\alpha)<+\infty$, and $\nu_{j,s,\alpha}=\delta_s$ for $\ell(\alpha)=+\infty$ and next, we define

$$
\nu_{j,s}=\delta_s \bigotimes \bigotimes_{\alpha \in D} \nu_{j,s,\alpha}.
$$

\Bin We have to prove that for any  $j\geq 1$, for $A \in \mathcal{A}$, the mapping

$$
s \mapsto \nu_{j,s}(A), \ \ (F1)
$$

\Bin is measurable. Let us show this in three steps. (i) We first show that for $A_\alpha \in \mathcal{S}_\alpha$ fixed, for $\alpha \in D$, for $j\geq 1$,

$$
s \mapsto \nu_{j,s,\alpha}(A_\alpha) \ \  (F2)
$$

\Bin is measurable. For $\ell(\alpha)=+\infty$, this mapping is the indicator function of $A_\alpha$. For $\ell(\alpha)<+\infty$, it is a constant function for $j>\ell(\alpha)$ or, for $1\leq j \leq \ell(\alpha)$, it is an elementary function associated to the subdivision of $S$ into the $C_{p,\ell(\alpha)}$'s. So the function in (F2) est measurable (in $s \in S$). (ii) Next, we know that the product $\sigma$-algebra $\mathcal{A}$ is generated by the class $\mathcal{C}$ of all cylinders of the form

$$
A=A_{\infty} \times \prod_{r=1}^{p} A_{\alpha_{r}} \times S^{\prime}_{\{\infty, \alpha_{1}, \cdots, \alpha_{p}\}},
$$ 

\noindent (where $p$ is finite, $A_{\infty} \in \mathcal{S}_{\infty}$, $A_{\alpha_p} \in \mathcal{S}_{\alpha_p}$, $1\leq r \leq p$). For such cylinders, we have, for $j\geq 1$ fixed,

$$
\nu_{j,s}(A)=1_{A_{\infty}}(s) \times \prod_{r=1}^{p} \nu_{j,s,\alpha}(A_{\alpha_r}),
$$

\Bin and hence, the mapping in (F1) is measurable for all $A \in \mathcal{C}$. (iii) But $\mathcal{C}$ is a $\pi$-system containing the full space $\Omega$. By \cite{ips-mestuto-ang}, Chapter 2, Exercise 4 in Doc 01-04), it is enough that show the class

$$
\{A \in \mathcal{A}, \ s\mapsto \nu_{j,s}(A) \ is \ measurable\}
$$  

\Bin is a $\lambda$-system containing $\mathbb{C}$. But this is quite direct and left as an exercise. So we may integrate $\nu_{j,s}(A)$ over $s \in S$ and define

$$
\mathcal{A} \ni A \mapsto \nu_j = \int_{S} \nu_{j,s}(A) \ d\mathbb{P}_{\infty}(s), 
$$

\Bin so that, just by the monotone convergence theorem, $\nu_j$ is a probability measure on $(\Omega, \mathcal{A})$ for each $j\geq 1$. Finally, we define

$$
\nu=\sum_{j\geq 1} \beta_j \nu_j.
$$

\bigskip \noindent We got the probability measure $\nu$ and the random variables $X_{\infty}$ and $X_{\alpha}$ we were searching and we have :\\

\begin{lemma} \label{fact7} $\nu$ is a probability measure and : \\

\noindent (i) $\nu X_{\infty}^{-1}=\mathbb{P}_{\infty}$.\\

\noindent (ii) For all $\alpha \in D$ with $\ell(\alpha)<+\infty$, $\nu X_{\alpha}^{-1}=\mathbb{P}_{\alpha}$.\\

\noindent (iii) For all $\alpha \in D$ with $\ell(\alpha)=+\infty$, $\nu X_{\alpha}^{-1}=\mathbb{P}_{\infty}$ (on $S_{\alpha})$. \ $\Diamond$\\
\end{lemma}

\Ni \textbf{Proof of Lemma \ref{fact7}}. The mapping $\nu$ is a probability measure as a non-negative linear combination of probability measure associated to constants which add up to one. Further, we have :\\

\noindent (i) We have for all $A_{\infty} \in \mathcal{S}_{\infty}$,  

\begin{eqnarray*}
\nu X_{\alpha}^{-1}(A_{\infty})&=&\nu\left(A_{\infty} \times S^{\prime}_{\{\infty\}}\right)\\
&=& \sum_{j\geq 1} \beta_j \int_S \left(1_{A_{\infty}}(s)\right) \ d\mathbb{P}_{\infty}(s)\\
&=& \sum_{j\geq 1} \beta_j \mathbb{P}_{\infty}(A_{\infty})\\
&=& \mathbb{P}_{\infty}(A_{\infty}).\\
\end{eqnarray*}

\noindent (ii) If $\ell(\alpha)<+\infty$, we have for all $A_{\alpha} \in \mathcal{S}_{\alpha}$, 

\begin{eqnarray*}
\nu X_{\alpha}^{-1}(A_{\alpha})&=&\nu\left(A_{\alpha} \times S^{\prime}_{\{\alpha\}}\right)\\
&=& \sum_{j\geq 1} \beta_j \int_S  \nu_{j,s,\alpha}(A_\alpha) \ d\mathbb{P}_{\infty}(s)\\
&=& \sum_{j>\ell(\alpha)} \beta_j \int_S  \nu_{j,s,\alpha}(A_\alpha) \ d\mathbb{P}_{\infty}(s)\\
&+& \sum_{j\leq \ell(\alpha)} \beta_j \int_S  \nu_{j,s,\alpha}(A_\alpha) \ d\mathbb{P}_{\infty}(s),\\
\end{eqnarray*}

\noindent which, when combined with the definition of $\nu_{j,s,\alpha}$ for $\ell(\alpha)<+\infty$, $\nu X^{-1}_{\alpha}(A_{\alpha})$ leads to 

\begin{eqnarray*}
\nu X_{\alpha}^{-1}(A_{\alpha})&=&
 \sum_{j>\ell(\alpha)} \beta_j \int_S  \mathbb{H}_{\alpha}(A_{\alpha}) \ d\mathbb{P}_{\infty}(s)\\
&+& \sum_{j\leq \ell(\alpha)} \beta_j \sum_{)\leq p\leq q(\ell(\alpha))} \int_{C_{p,\ell(\alpha)}}  \nu_{j,s,\alpha}(A_\alpha) \ d\mathbb{P}_{\infty}(s)\\
&=& (1-\beta^{\ast}_{\ell(\alpha)})  \mathbb{H}_{\alpha}(A_{\alpha}) \\
&+& \sum_{j\leq \ell(\alpha)} \beta_j \sum_{1\leq p\leq q(\ell(\alpha))} \int_{C_{p,\ell(\alpha)}} \mathbb{P}_{\alpha}(A_{\alpha}/C_{p,\ell(\alpha)}) \ d\mathbb{P}_{\infty}(s)\\
&=&(1-\beta^{\ast}_{\ell(\alpha)})  \mathbb{H}_{\alpha}(A_{\alpha}) + \beta^{\ast}_{\ell(\alpha)} \sum_{0\leq p \leq q(\ell(\alpha))} \mathbb{P}_{\alpha}(A_{\alpha}/C_{p,\ell(\alpha)}) \mathbb{P}_{\infty}(C_{p,\ell(\alpha)}).
\end{eqnarray*}

\bigskip \noindent By comparing the later line and Formula \eqref{wcib_defHalphaB}, we get that $\nu X_{\alpha}^{-1}=\mathbb{P}_{\alpha}$. Finally for $\ell(\alpha)=+\infty$, we have

\begin{eqnarray*}
\nu X_{\alpha}^{-1}(A_{\alpha})&=&\nu\left(A_{\alpha} \times S^{\prime}_{\{\alpha\}}\right)\\
&=& \sum_{j\geq 1} \beta_j \int_S  \nu_{j,s,\alpha}(A_\alpha) \ d\mathbb{P}_{\infty}(s)\\
&=& \sum_{j\geq 1} \beta_j \int_S  1_{A_\alpha}(s) \ d\mathbb{P}_{\infty}(s)\\
&=& \sum_{j\geq 1} \beta_j \mathbb{P}_{\infty}(A_\alpha)\\
&=&\mathbb{P}_{\infty}(A_\alpha).\\
\end{eqnarray*}

\bigskip \noindent We are ready to conclude. Let us defined for $k\geq 1$,

$$
A\equiv \bigcap_{k\geq 1} A_k \equiv \bigcap_{k\geq 1} \bigcup_{\ell(\alpha)\geq k} (d(X_\alpha,X_{\infty})> \Delta_k).
$$

\bigskip \noindent For a fixed $\alpha \in D$, $\ell(\alpha)\geq k$, we define $G_{k,\alpha}=(d(X_\alpha,X_{\infty})>\Delta_k)$. We do not know whether $A$ or the $A_k$'s are measurable or not. But, by Proposition \ref{wcib_07_countable_Product}, there exists a measurable set $A^\ast_k$ which includes $A_k$ and depends only on a countable set $D_{0,k}$ of indices, so that for $D_{0}$ being the union of the $D_{0,k}$'s to which we add $\{\infty\}$ if needed , we have

$$
A_k \subset A^\ast_k=A^\ast_{k,{D_0}} \times S^{\prime}_{D_0}, \ k\geq 1 \ \ (CC1)
$$

\bigskip \noindent and by applying $\Pi_{D_0}$ to that formula, we have

$$
\Pi_{D_0}(A_k) \subset A^\ast_{k,{D_0}}. \ \ \ (CC2)
$$

\Bin We also have by taking complements in \textit{(CC1)}, 

$$
\biggr(\bigcap_{\ell(\alpha)\geq k} \left(A^\ast_{k,{D_0}}\right)^c\biggr) \times S^{\prime}_{D_0} \subset A_k^c. \ \ \ (CC3)
$$

\bigskip \noindent But we have, as a general rule

$$
A_k \subset \pi_{D_0}(A_k) \times S^{\prime}_{D_0}. \ \ \ (CC4)
$$

\Bin The combination of Formulas (CC1)-(CC4) leads to

$$
\left(A^\ast_{k,{D_0}}\right)^c \times S^{\prime}_{D_0} \subset \pi_{D_0}(A_k^c) \times S^{\prime}_{D_0} \subset A_k^c.
$$

\bigskip \noindent Besides, it is clear that 

$$
\pi_{D_0}(A_k^c) \times S^{\prime}_{D_0}=\biggr(\bigcap_{\ell(\alpha)\geq k, \alpha \in D_0} G_{k,\alpha}^c\biggr) \times S^{\prime}_{D_0}.
$$

\bigskip \noindent The combination of the two later facts yields

$$
A_k \subset \biggr(\bigcup_{\alpha \in D_0} G_{k,\alpha}\biggr) \times S^{\prime}_{D_0} \subset A_k^\ast.
$$

\bigskip \noindent Let us study

$$
A^{\ast\ast}_k \equiv \bigcup_{\ell(\alpha)\geq k, \ \alpha \in D_0} (d(X_\alpha,X_{\infty}) > \Delta_k).\ \ \ (F1)
$$

\noindent Let us fix $k\geq 1$ and $\alpha \in D$ and set $G_{k,\alpha}=(d(X_\alpha,X_\infty)>\Delta_k)$. We have, for any $h\geq \ell(\alpha)\geq k$  

\begin{eqnarray*}
\nu_j(G_{k,\alpha})&=&\int_S \nu_{j,s}(G_{k,\alpha}) \ d\mathbb{P}_{\infty}(s)\\
&=&\sum_{0\leq p \leq q(\ell(\alpha))} \int_{C_{p,\ell(\alpha)}} \nu_{j,s}(G_{k,\alpha}) \ d\mathbb{P}_{\infty}(s)\\
&=&\sum_{0\leq p \leq q(\ell(\alpha))}  \int_{C_{p,\ell(\alpha)}} \nu_{j,s}((X_{\infty}\in C_{0,h}) \cap G_{k,\alpha}) \ d\mathbb{P}_{\infty}(s)\\
&+&\sum_{0\leq p \leq q(\ell(\alpha))} \sum_{1\leq u \leq q(h)} \int_{C_{p,\ell(\alpha)}} \nu_{j,s}((X_{\infty}\in C_{u,h}) \cap G_{k,\alpha}) \ d\mathbb{P}_{\infty}(s)\\
&=:&T_1+T_2.
\end{eqnarray*}

\Bin We have

\begin{eqnarray*}
T_1&\leq &\sum_{0\leq p \leq q(\ell(\alpha))}  \int_{C_{p,\ell(\alpha)}} \nu_{j,s}(X_{\infty}\in C_{0,h})  \ d\mathbb{P}_{\infty}(s)\\
&=&\sum_{0\leq p \leq q(\ell(\alpha))}  \int_{C_{p,\ell(\alpha)}} \delta_s(C_{0,h})  \ d\mathbb{P}_{\infty}(s)\\
&=&\sum_{0\leq p \leq q(\ell(\alpha))}  \int_{C_{p,\ell(\alpha)}} 1_{C_{0,h}}(s) \ d\mathbb{P}_{\infty}(s)\\
&=&\sum_{0\leq p \leq q(\ell(\alpha))}  \mathbb{P}_{\infty}(C_{p,\ell(\alpha)} \cap C_{0,h})\\
&=&  \mathbb{P}_{\infty}(C_{0,h})\\
&\leq & \varepsilon_{h}.
\end{eqnarray*}

\Bin Next, we can use Fubini's theorem to see that

\begin{eqnarray*}
T_2&=&\sum_{0\leq p \leq q(\ell(\alpha))} \sum_{1\leq u \leq q(h)} \int_{C_{p,\ell(\alpha)}} 
\int_{S_{\infty}} 1_{C_{u,h}}(x_{\infty}) \biggr(\int_{S_{\alpha}}  \nu_{j,s,\alpha} (B^c(x_{\infty},\Delta_k)) d\delta_s(x_\infty)\biggr) d\mathbb{P}_{\infty}(s)\\
&=&\sum_{0\leq p \leq q(\ell(\alpha))} \sum_{1\leq u \leq q(h)} \int_{C_{p,\ell(\alpha)}} 1_{C_{u,h}(s)} \nu_{j,s,\alpha}(B^c(s,\Delta_k)) d\mathbb{P}_{\infty}(s).
\end{eqnarray*}

\Bin By construction,  each $C_{u,h}$ is in some $C_{p,\ell(\alpha)}$. Hence, we may define for $0\leq p \leq q(\ell(\alpha))$,

$$
I(p,h)=\{u, \ u \in \{1,\cdots,q(h)\} \ and \ C_{u,h}\subset C_{p,\ell(\alpha)}\}.
$$ 

\Bin So for $u \notin I(p,h)$, $C_{u,h} \cap C_{p,\ell(\alpha)} = \emptyset$. Thus, for

$$
B(p,h)=\sum_{u \in I(p,h)} C_{u,h},
$$

\Bin the last equation reduces to 

\begin{eqnarray*}
T_2=\sum_{0\leq p \leq q(\ell(\alpha))}  \int_{B(p,h)}  \nu_{j,s,\alpha}(B^c(s,\Delta_k)) d\mathbb{P}_{\infty}(s).
\end{eqnarray*}

\Bin Let

$$
T_2=\sum_{0\leq p \leq q(\ell(\alpha))} T(2,p) \ with \ T(2,p)=\int_{B(p,h)}  \nu_{j,s,\alpha}(B^c(s,\Delta_k)) d\mathbb{P}_{\infty}(s).
$$

\Bin We have

\begin{eqnarray*}
T(2,p)&=& \int_{B(p,h)} \nu_{j,s,\alpha}(B^c(s,\Delta_k) \cap C_{p,\ell(\alpha)}) d\mathbb{P}_{\infty}(s)\\
&+&\int_{C_{p,\ell(\alpha)}}  \nu_{j,s,\alpha}(B^c(s,\Delta_k) \cap C^{c}_{p,\ell(\alpha)}) d\mathbb{P}_{\infty}(s)\\
&=&T(21,p) + T(22,p).
\end{eqnarray*}

\Bin Now, in what follows, we may do the computations of the probability space $(S^2, \mathcal{A}^{\otimes 2}, \mathbb{P}_{\infty} \times \nu_{j,s,\alpha})$ since the measurable spaces $S_\infty$ and $S_{\alpha}$ are identical. So, by Tonelli's theorem, se have

\begin{eqnarray*}
T(21,p)&=& \int_{S} 1_{B(p,h)}(s_1) \int_{S}  1_{C_{p,\ell(\alpha)}}(s_2)  1_{B^c(s_1,\Delta_k)}(s_2) \ d\mathbb{P}_{\infty}(s_1) \ d\nu_{j,s_1,\alpha}(s_2)\\
&=& \int_{S\times S} \biggr(1_{B(p,h)}(s_1) 1_{B^c(s_1,\Delta_k)}(s_2) \ d\mathbb{P}_{\infty}(s_1)\biggr) 1_{C_{p,\ell(\alpha)}}(s_2) \ d\nu_{j,s_1,\alpha}(s_2),
\end{eqnarray*}

\Bin in which the expression between the big parentheses is zero since the diameter of $C_{p,\ell(\alpha)}$ is less or equal to 
$\Delta_{\ell(\alpha))}\leq \Delta_{k}$ and  for $(s_1,s_2) \in C^2_{p,\ell(\alpha)}$. So $T(21,p)=0$.

\Bin Next, by  Tonelli's theorem again,
$$
T(22,p)= \int_{S} 1_{B(p,h)}(s_1) \biggr(\int_{S}  1_{C^c_{p,\ell(\alpha)}}(s_2)  1_{B^c(s_1,\Delta_k)}(s_2) d\nu_{j,s_1,\alpha}(s_2) \biggr)\ d\mathbb{P}_{\infty}(s_1).\\
$$

\Bin Either $\ell(\alpha)=+\infty$ and in that case, $\nu_{j,s_1,\alpha}=\delta_{s_1}$ and hence, the integral between the big parentheses, is

$$
1_{C^{c}_{p,\ell(\alpha)}}(s_1)  1_{B^c(s_1,\Delta_k)}(s_1)=0,
$$

\Bin and hence $T(22,p)=0$. Or $\ell(\alpha)<\infty$ and in that case, $\nu_{j,s_1,\alpha}(\circ)=\mathbb{P}_{\alpha}(\circ/C_{p,\ell(\alpha)})$ (we recall that
$B(p,h)\subset C_{p,\ell(\alpha)})$ and hence, the integral between the big parentheses, is

$$
\mathbb{P}_{\alpha}(\{B^c(s_1,\Delta_k) \cap C^c_{p,\ell(\alpha)}\} / C_{p,\ell(\alpha)})=0. 
$$
 
\Bin We conclude that $T_2=0$ and thus $\nu_j(G_{k,\alpha}) \leq \varepsilon_{h}$ and finally, for any $\alpha \in D$, for any $k\geq 1$, for any $h\geq \ell(\alpha)\geq k$,

$$
\nu(G_{k,\alpha})\leq \varepsilon_{h}.  (CC)
$$

\Bin We need a little extra-work to do before concluding. For $k\geq 1$ fixed, we denote the countable set  $\{\alpha \in D_0, \ \ell(\alpha)\geq k\}$ by 
$\{\alpha_1, \alpha_2, \cdots \}$.\\

\Bin We are going to use the construction that let Formula (CC) in an induction reasoning. For $r=1$ we choose $h(1)$ such that

$$
\nu(G_{\alpha_1,k})\leq \varepsilon_{h(1)}.
$$

\Bin Next for $r=2$, we choose $h(2)> max(h(1), \ell(\alpha_2))$ to get

$$
\nu(G_{\alpha_2,k})\leq \varepsilon_{h(2)}
$$

\Bin By induction, we get $k<h(1)<\cdots<h(2)<\cdots$ such that for any $r\geq 1$,

$$
\nu(G_{\alpha_r,k})\leq \varepsilon_{h(r)}.
$$

\Bin We conclude that

$$
\nu\biggr(\bigcup_{\ell(\alpha)\geq k, \ \alpha \in D_0} (d(X_\alpha,X_{\infty}) > \Delta_k)\biggr)\leq \sum_{r\geq 1} \varepsilon_{h(r)}, 
$$
 
\Bin i.e.,

$$
\nu\biggr(A^{\star \star}_k\biggr) \leq \sum_{j\geq k} \varepsilon_{j}. 
$$

\Bin Since the series $\sum \varepsilon_j$ converges, we get that

$$
\nu\left( A^{\star \star}\right)=0, 
$$

\Bin where

$$
A^{\star \star}=\bigcap_{k\geq 1} A^{\star \star}_k.
$$

\Bin By definition, the set $(X_{\alpha} \rightarrow X_{\infty})$ includes $A^c$, and we have

$$
A \subset A^{\star \star}.
$$

\Bin So the exterior measure of $(X_{\alpha} \rightarrow X_{\infty})^c$ is given by

$$
\nu^{\star}(X_{\alpha} \nrightarrow X_{\infty})=0.
$$

\Bin In a last move, we may extend $(\Omega_D,\mathcal{S}_D,\nu)$ to a complete probability measure
$(\Omega_D, \tilde{\mathcal{S}}_D,\tilde{\nu})$ (See \cite{ips-mestuto-ang}, Chapter 5, Doc 04-02, Exercise 8). We will have

$$
\tilde{\nu}\biggr( \limsup_{k\rightarrow +\infty} \bigcup_{\ell(\alpha)\geq k, \alpha \in D}(d(X_{\alpha}, X_{\infty}) <\Delta_k)\biggr)=0.
$$

\bigskip \noindent Since the mappings $X_{\alpha}$, $\alpha \in \{\infty\} \cup D$, do not change and all the previous laws concern $\mathcal{S}_D$-measurable mappings or sets (for the handling of which, $\nu$ and $\tilde{\nu}$ are equivalent), we conclude that

$$
d(X_{\alpha}, X_{\infty}) \rightarrow_{D} 0, \ a.s..
$$

\bigskip \noindent The proof is over. $\blacksquare$\\

\noindent \textbf{Remarks}.\\

\noindent \textbf{(R1)} For interested readers who want to compare with the original proof in \cite{wichura1970}, the given proof is more direct in the sense that
we did not use the step where $D$ is countable and $S$ is a finite metric space in Section 2. Neither we did require the existence of the mapping $k(\gamma)$ of $\gamma \in D$ in Formula (e) therein. In our approach, such a sequence came to birth itself.\\

\noindent \textbf{(R2)} We recall a very simple proof of the Skorohod theorem for $\mathbb{R}^d$ with $d=1$. It would be interesting to either find out if simple proofs are available or can be done for $d\geq 2$.\\

\noindent \textbf{(R3)} The present contribution is the result of a Msc Dissertation at University Gaston Berger (SENEGAL).

\end{document}